\documentclass[12pt]{amsart}
\usepackage[centertags]{amsmath}
\usepackage{amsfonts}
\usepackage{amsthm}
\usepackage{newlfont}
\usepackage{amscd}
\usepackage{amsgen}
\usepackage{amssymb}
\newlength{\defbaselineskip} 
\theoremstyle{plain} \newtheorem{thm}{Theorem}[section]
\newtheorem{pro}[thm]{Problem} \newtheorem{lem}[thm]{Lemma}
\newtheorem{prop}[thm]{Proposition}
\theoremstyle{definition}
\newtheorem{defi}[thm]{Definition}
\newtheorem{rem}[thm]{Remark}

 \numberwithin{equation}{section}
\numberwithin{equation}{section}
\DeclareMathOperator{\Pic}{Pic}
 \DeclareMathOperator{\Spec}{Spec}
 \DeclareMathOperator{\Def}{Def}
\DeclareMathOperator{\sing}{sing}
 \DeclareMathOperator{\coker}{coker}
\DeclareMathOperator{\im}{im}
\begin{document}

\title{Primitive contractions of Calabi-Yau threefolds I}
\author{Grzegorz Kapustka}
\author{Micha\l\ Kapustka}
\thanks{The authors  are scholars of the project which is
co-financed from the European Social Fund and national budget in
the frame of The Integrated Regional Operational Programme.}
\begin{abstract}
We construct examples of primitive contractions of Calabi--Yau
threefolds with exceptional locus being $ \mathbb{P}^1 \times
\mathbb{P}^1$, $\mathbb{P}^2$, and smooth del Pezzo surfaces of
degrees $\leq 5$. We describe the images of these primitive
contractions and find their smoothing families. In particular, we
give a method to compute the Hodge numbers of a generic fiber of
the smoothing family of each Calabi--Yau threefold with one
isolated singularity obtained after a primitive contraction of
type II. As an application, we get examples of natural conifold
transitions between some families of Calabi--Yau threefolds.
\end{abstract}
\maketitle
\section{Introduction}
A primitive Calabi--Yau contraction (cf.~\cite{Wilson2}) is a
birational morphism between Calabi--Yau varieties that does not
factorize in the category of algebraic varieties. We say that a
primitive contraction is of type II, if its exceptional locus is
an irreducible surface that is mapped to a point. The dualizing
sheaf of the exceptional surface is then anti-ample (i.e. the
surface is a generalized del Pezzo surface). Moreover, from
\cite[thm.~5.2]{Gross1}, we know that this exceptional locus is
either a normal (Gorenstein) rational del Pezzo surface or a
non-normal Gorenstein surface with $\omega_E^{-1}$ ample of degree
$7$.

This paper grew up from an attempt to understand which del Pezzo
surfaces do occur as the exceptional locus of a primitive
contraction. We give examples of primitive contractions with
exceptional locus being $\mathbb{P}^2$,
$\mathbb{P}^1\times\mathbb{P}^1$ or a smooth del Pezzo surface of
degree $\leq 5$. To obtain these examples, we need to find a smooth
Calabi--Yau threefold $X$ containing the chosen del Pezzo surface
$D$ in such a way that each divisor on $X$ restricts to a multiple
of the canonical divisor on $D$. Two types of constructions are
given.

In the first, we consider the resolutions of double coverings of Fano
varieties with index $\geq 2$, branched along the sum of two smooth
surfaces $G'\in |qK_F|$ and $D'\in |(2-q)K_F|$, where $q$ is chosen
such that $D'$ is a smooth del Pezzo surface. The strict transform
$D$ of $D'$ will be an exceptional locus of a primitive contraction
of type II. In this way, we obtain primitive contractions of all
considered del Pezzo surfaces.

For del Pezzo surfaces of degree $3$, $4$, and $5$, we find
another construction (constructions 3 and 4). The Calabi--Yau
threefolds are then obtained as small resolutions of nodal
threefolds (a quintic and a complete intersection of two cubics,
respectively) containing the anti-canonical models of the
considered del Pezzo surfaces.

The examples of primitive contractions of del Pezzo surfaces of
degrees $6$, $7$, of $\mathbb{P}^2$ blown up in one point, and of
singular del Pezzo surface are more difficult to find. This
problem is treated in \cite{K}.

The second aim of this paper is to describe the smooth fibers of
smoothing families of singular Calabi--Yau threefolds obtained by
constructions contained in this paper. More precisely, the image
of a primitive contraction of type II is a Calabi--Yau threefold
with one isolated rational Gorenstein (i.e. canonical index $1$),
$\mathbb{Q}$-factorial singularity. By the results of Gross
(\cite[Thm.~5.8]{Gross1}), for most such Calabi--Yau threefolds,
there exists a flat family over the unit disc such that the
generic fiber is smooth and the central fiber is this manifold.
Such a family is called a smoothing of the singular Calabi--Yau
threefold. We describe smooth fibers of these families (they are
projective Calabi--Yau threefolds if the central threefold is
projective).

In section \ref{smoothing}, we describe the miniversal Kuranishi
space of each singular Calabi--Yau threefold obtained as the image
of a primitive contraction of type II. This gives us a method to
compute the Hodge numbers of smoothings of these singular
varieties.

We give moreover an exact description, with equations in projective
spaces, of the singular Calabi--Yau threefolds considered in our
constructions. To do this, we first find the linear system that
gives the primitive contraction and then describe the image of the
morphism associated to this linear system.

In this paper, we also study our examples in the context of
conifold transition. The idea of conifold transitions goes back to
Miles Reid (see \cite{Reid}) and was generalized in \cite{Gross2}
in the following way. We say that the moduli spaces of two smooth
Calabi--Yau 3-folds $X$ and $Y$ are joined by a conifold
transition, if there exists a Calabi--Yau threefold $Z$ and two
birational contractions $Z\rightarrow X'$, $Z\rightarrow Y'$ such
that $X'$ (resp.~$Y'$) has a smoothing inside the moduli space of
$X$ (resp.~$Y$) (i.e. there exists a flat family
$\mathcal{X}\rightarrow \Delta$ such that $\mathcal{X}_0 \simeq
X'$ and $\mathcal{X}_t \in X$ for $t\neq 0$). It is conjectured
that all families of Calabi--Yau threefolds can be joined by a
sequence of conifold transitions.

All obtained results are summarized in Table \ref{table intr.},
where the following notation is used. The symbols
$X_{d_1,d_2,\dotsc}$ denote complete intersections of indicated
degrees in the indicated manifold, whereas $X^2$ means that $X$ is
embedded by the double Veronese embedding. Moreover,
$X\xrightarrow{2:1}F_i$ denotes a double cover of a Fano threefold
$F_i$ of index $2$ branched over $-2K_{F_i}$, $Q_3$ is a
three-dimensional quadric, and $n\times n$ Pffafian stands for the
variety defined by the $n-1 \times n-1$ Pffafians of an
appropriate $n\times n$ skew-symmetric matrix.

 \begin{table}[h]\label{table intr.}
 \caption{collected results}
             \begin{center}

            \renewcommand*{\arraystretch}{1}
            \begin{tabular}{|c|c|c|c|c|c|}\hline
            $\deg D$&$X'$&$C=\sing X'$&$\chi(\mathcal{Y}_t)$&$Y'$& thm.                    \\ \hline
                                                                         \hline

            1&$X'\xrightarrow{2:1}F_1$                       &$C_{3,1}\subset F_1$&-204      &$Y_6\subset{P}(1,1,1,1,2)$           &  4.2 \\ \hline
            2&$X'\xrightarrow{2:1}F_2$                       &$C_{3,1}\subset F_2$&-156      &$Y_{3,4}\subset \mathbb{P}(1,1,1,1,1,2)$&   4.2\\ \hline
            3&$X'\xrightarrow{2:1}F_3$                       &$C_{3,1}\subset F_3$&-144      &$Y_{3,3}\subset \mathbb{P}^5$        &   4.2\\ \hline
            4&$X'\xrightarrow{2:1}F_4$                       &$C_{3,1}\subset F_4$&-144      &$Y_{2,2,3}\subset \mathbb{P}^6$      &  4.2\\ \hline
            5&$X'\xrightarrow{2:1}F_5$                       &$C_{3,1}\subset F_5$&-150      &$Y_{3,1,1}\subset G(2,5)$                     &  4.2 \\ \hline
            8&$X'\xrightarrow{2:1}\mathbb{P}^3$              &$C_{6,2}\subset \mathbb{P}^3$&-204      &$Y_6\subset\mathbb{P}(1,1,1,1,2)$  &  4.3 \\ \hline
            9&$X'\xrightarrow{2:1}\mathbb{P}^3$              &$C_{7,1}\subset \mathbb{P}^3$& -     &$\deg (63)\subset\mathbb{P}^{20}$      &  4.3 \\ \hline
            3&$X'\xrightarrow{2:1}\mathbb{P}^3$              &$C_{5,3}\subset \mathbb{P}^3$&-200      &$Y_5\subset\mathbb{P}^4$             &  4.3 \\ \hline
            4&$X'\xrightarrow{2:1} Q_3$                      &$C_{2,2,4}\subset \mathbb{P}^4$&-176      &$Y_{2,4}\subset\mathbb{P}^5$         &  4.4 \\ \hline
            8&$X'\xrightarrow{2:1} Q_3$                      &$C_{1,2,5}\subset \mathbb{P}^4$&-200      &$Y^2_5\subset\mathbb{P}^{14}$        &  4.4\\ \hline
            3&$X'_5\subset\mathbb{P}^4$                      &24 ODP&-176      &$Y_{2,4}\subset\mathbb{P}^5$                              &  5.2 \\ \hline
            4&$X'_5\subset\mathbb{P}^4$                      &36 ODP&-144      &$Y_{3,3}\subset\mathbb{P}^5$                              &  5.1 \\ \hline
            5&$X'_{3{,}3}\subset\mathbb{P}^5$                &28 ODP&-98      &$7\times7 \operatorname{Pffafian}\subset \mathbb{P}^6$             & 5.3\\ \hline
            3&$X'_{3{,}3}\subset\mathbb{P}^5$                &12 ODP&-144      &$Y_{2,2,3}\subset \mathbb{P}^6$                          &  5.4 \\ \hline
            4&$X'_{3{,}3}\subset\mathbb{P}^5$                &20 ODP&-120      &$5\times5 \operatorname{Pffafian}\subset\mathbb{P}^6$            &  5.4 \\ \hline

                       \end{tabular}

            \end{center}

            \end{table}
\section*{Acknowledgements}
We would like to thank S.~Cynk for his constant help. Some ideas
from this paper grew up during our stay at the EAGER training site
at Warwick University in spring $2004$, we would like to express
our gratitude to M.~Reid for mathematical inspiration. We thank
M.~Gross, A.~Langer, V.~Nikulin, P.~Pragacz, J.~Wi\'{s}niewski for
answering questions. We acknowledge S.~Kleiman for help in the
redaction.
\section{Primitive contractions}
By a Calabi--Yau threefold $X$ we mean a complex projective
threefold with canonical singularities such that the canonical
divisor $K_X=0$ and $h^1(\mathcal{O}_X)=0$. The aim of this paper
is to study primitive contractions. Let us first recall some
definitions from \cite[p.~566]{Wilson2}.

            \begin{defi}

            Let $X$ be a smooth Calabi--Yau threefold. We say that a birational morphism $\phi\colon X\rightarrow Y$
            is a \emph{primitive contraction}, if $Y$ is normal and one of
            the following equivalent conditions holds:
            \begin{itemize}
            \item[a)] $\dim \Pic_{\mathbb{R}}(Y)=\dim
            \Pic_{\mathbb{R}}(X)-1$;

            \item[b)] $\phi$ does not factorize in the category of normal
            algebraic varieties.
            \end{itemize}
            \end{defi}

            \begin{defi}

            We say that a primitive contraction is
            \begin{itemize}
            \item[-] of \emph{type I}, if it
            contracts only finitely many curves;
            \item[-] of \emph{type II}, if it
            contracts an irreducible surface down to a point;
            \item[-] of
            \emph{type III}, if it contracts a surface down to a curve.
            \end{itemize}
            \end{defi}
 Our aim is to study type II contractions. First, summarizing the
 results from \cite{Gross1} and \cite{Wilson2} we obtain necessary conditions for the exceptional loci of such contractions.
\begin{prop}\label{GrossWilson} The exceptional divisor of a primitive contraction of
type II is a generalized del Pezzo surface (i.e. an irreducible
Gorenstein surface with anti-ample dualizing sheaf). Moreover, the
exceptional locus $E$ is either a rational Gorenstein del Pezzo
surface, or a non-normal surface with $\omega_E^2=7$ (that is not
a cone), or a cone over an elliptic curve such that
$\omega_E^2\leq3$.
\end{prop}
\begin{proof} To prove that the exceptional locus is a generalized
del Pezzo surface observe that $\pi\colon X\rightarrow Y$ is a
resolution of singularities. Since $\pi$ is indecomposable, we
conclude that $\pi$ is a blow-up or an $\alpha$-blow-up of the
singular point (to see this, we can argue as in the proof of
\cite[lemma~2.3]{Reid can}). Now from proposition 2.13 \cite{Reid
can} we obtain that the exceptional locus is a generalized del
Pezzo surface. The second part follows from
\cite[thm.~5.2]{Gross1}.
\end{proof}
This paper grew-up from an attempt to understand the following
problem.
\begin{pro}  Which smooth del Pezzo surfaces can be
contracted by a primitive contraction of type II.
\end{pro}
 To prove that a given contraction is primitive we will use the
following characterization.
\newpage
\begin{lem}\label{pierwsze}
A del Pezzo surface $D$ can be contracted by a primitive
contraction if and only if $D$ can be embedded in a smooth
Calabi--Yau manifold $X$ in such a way that for each Cartier
divisor $E\subset X$ we have $E|_{D}= aK_{D}$, for some rational
number $a$.

\end{lem}
\begin{proof} Let $H$ be an ample divisor on $X$.
            We claim that the $\mathbb{Q}$-divisor $aH+D$, where
            $a=\frac{(-K_D)^2}{H(-K_D)}$ (i.e.
            $aH|_{D}= -K_D$), is big and \emph{nef}. For this it is enough to observe that
            if $C$ is a curve contained in $D$ then $C.(aH+D)=0$. From the basepoint-free theorem
            we find an $n$ such
that $|n(aH+D)|$ gives a birational morphism. Next, we shall
follow the proof of \cite[thm.~1.2]{Wilson3} to prove that we can
find an integer $m$ such that the image of
$\varphi_{|mn(aH+D)|}=:\varphi_m $ is a normal variety such that
the morphism is an isomorphism outside from $D$. Let
$\varphi\colon X\xrightarrow{\gamma} Y\xrightarrow{\psi} Z$ be the
Stein factorization such that $\gamma$ has connected fibers and
$\psi$ is finite. The pull-back $\psi^*(H_1)=H$ of the hyperplane
section of $Z$ is an ample divisor on $Y$ (see \cite[exercise~5.7
d]{Har}). We choose an integer $r$ such that $rH=\psi^*(rH_1)$ is
very ample. Since the Stein factorization of the morphism
$\varphi_r$ factorizes through $\varphi_1$ the morphism
$\varphi_r$ is birational and has normal image.
\end{proof}
            \section{Smoothings}\label{smoothing}
Most of the Calabi--Yau threefolds obtained after a primitive
contraction of type II are smoothable.
M. Gross proved (see \cite[Thm.~5.8]{Gross1})
             that if $X\rightarrow Y$ be a primitive type
            II contraction with exceptional divisor $E$. Then $Y$ is smoothable unless
$E\simeq \mathbb{P}^2$ or $E\simeq \mathbb{F}_1$.

Our aim is to compute the Hodge numbers of the obtained smooth varieties.

\begin{prop}\label{h11}  Let $X\rightarrow Y$ be a primitive
contraction of type II, such that $Y$ is smoothable. If
$\mathcal{Y}\rightarrow \Delta$ is a $1$-parameter smoothing, then
$$ h^{1,1}(X)-1=h^{1,1}(\mathcal{Y}_t).$$
\end{prop}

\begin{proof}(cf.~\cite[prop.~6.1]{namikawa}) Let $S$ be the universal
Kuranishi space of $Y$ and $\mathcal{Z}\rightarrow S$ be the
universal family. It follows from the Artin approximation theorem
(see \cite[1.6]{Artin}) that we can find a flat projective
morphism of algebraic varieties
$\overline{\pi}:\overline{\mathcal{Z}}\rightarrow \overline{S}$
such that $S$ in an open subset of $\overline{S}$ and
$\overline{\pi}|_S=\pi$. From \cite[12.1.10~and~12.1.9]{K-M}, we
can assume that $\overline{\mathcal{Z}}$ is $\mathbb{Q}$-factorial
and that $\overline{\pi}$ has a section. This implies that every
irreducible component of the relative Picard scheme
$Pic_{\overline{\mathcal{Z}}|\overline{S}}$ is proper over
$\overline{S}$ (see \cite[p.~68]{Kl} and \cite[12.1.8]{K-M}).
Since $H^1(\mathcal{O}_X)=H^2(\mathcal{O}_X)=0$, from the general
theory of Picard schemes (see \cite[prop.~5.19]{Kl}), we obtain
that each irreducible component of
$Pic_{\overline{\mathcal{Z}}|\overline{S}}$ is \'{e}tale over
$\overline{S}$ (after shrinking $\overline{S}$ to a suitable
Zariski open subset). Thus, over $S$ the Picard scheme
$Pic_{\overline{\mathcal{Z}}|\overline{S}}\times_{\overline{S}}S=Pic_{\mathcal{Z}|S}$
is a disjoint union of countably many copies of $S$. The natural
structure group of $Pic_{\mathcal{Z}|S}$ induces the group
structure on the fibers. It follows that the ranks of the Picard
groups satisfy
\begin{equation*}
\rho(\mathcal{Y}_t)=\rho(\mathcal{Y}_0)=\rho(X)-1.\qedhere
\end{equation*}
\end{proof}
\begin{rem} The above proposition follows also from \cite[prop.~3.1]{Gross2} and
\cite[12.2.1.4.2]{K-M}, since $Y$, as the image of a primitive
contraction of type II, has rational $\mathbb{Q}$-factorial
singularities. Note that this proposition does not work for
primitive contractions of type III. Indeed, since we cannot use
\cite[12.1.9]{K-M}, we cannot prove the properness of the
irreducible components of
$Pic_{\overline{\mathcal{Z}}|\overline{S}}$. See \cite[ch.~
3]{Gross2} for an explicit example.
\end{rem}
We first compute $h^{1,2}(\mathcal{Y}_t)-h^{1,2}(X)$ in the case
when $Y$ has non-hypersurface singularities.
\begin{thm}\label{hodge}  Let $\pi\colon X\rightarrow Y$ be a
primitive contraction of type II, contracting a smooth del Pezzo
surface of degree $r$ such that $5\leq r \leq 8$. If
$\mathcal{Y}\rightarrow \Delta$ is a $1$- parameter smoothing,
then for $t\neq 0$ we obtain:
\begin{enumerate}
\item if $r=5$, then $h^{1,2}(\mathcal{Y}_t)= h^{1,2}(X)+4$;

\item if $r=6$, then $h^{1,2}(\mathcal{Y}_t)= h^{1,2}(X)+1$
 or $h^{1,2}(\mathcal{Y}_t)= h^{1,2}(X)+2$;

\item if $r=7$, then $h^{1,2}(\mathcal{Y}_t)= h^{1,2}(X)+1$;

\item if $r=8$, then $h^{1,2}(\mathcal{Y}_t)= h^{1,2}(X)+1$.
\end{enumerate}
\end{thm}

\begin{proof} Recall that, a miniversal Kuranishi space exists for
compact complex spaces. For a Calabi--Yau threefold $X$ with
canonical singularities this space is universal since
$\mbox{Hom}(\Omega ^1 _X,\mathcal{O}_X)=0$ (see \cite[cor.~
8.6]{Kaw}). We denote by $\Def(X)$ the analytic germ of the
Kuranishi space of $X$ or its suitable analytic representative. It
follows from the openness of universality (see \cite[Satz~
7.1(1)]{Bin}), that the possible values of
$h^{1,2}(\mathcal{Y}_t)$ will be equal to the dimensions of the
irreducible components of $\Def(Y) $.

To compute these dimensions, we consider $T^1$ and $T^1_{loc}$ the
tangent spaces to $\Def(Y)$ and $\Def(Y',P)$ respectively, where
$\Def(Y',P)$ denotes the miniversal deformation space of the germ
$(Y',P)$ of the singularity obtained after the primitive
contraction. From \cite[thm.~2]{Schl} we have
$$T^1 \simeq H^1 (Y-{P}, \Theta_Y) \simeq H^1(X-E,\Theta_X)$$ and
$$T^1_{loc} \simeq H^1(Y'-P,\Theta_{Y'})\simeq
H^1(X'-E,\Theta_{X'}).$$ Here $(X',E)\rightarrow (Y',P)$ is the
given minimal resolution.

We have the following exact local cohomology sequences

\[H^1(X' ,\Theta_{X'})\rightarrow H^1(X'-E ,\Theta_{X'})
\xrightarrow{f'} H^2_E (X', \Theta_{X'}) \rightarrow
H^2(X',\Theta_{X'})\]
$$\ \ \ \uparrow \ \ \ \ \ \ \ \ \ \ \ \ \ \ \ \ \ \ \ \ \ \ \ ||$$
$$H^1(X ,\Theta_{X})\ \rightarrow \ H^1(X-E ,\Theta_{X}) \ \xrightarrow{f}
\ H^2_E (X, \Theta_{X})  \xrightarrow{\varphi} \
H^2(X,\Theta_{X}).$$

\medskip
It follows from \cite[lem.~4.5]{Gross1} that $\im f' =\im f=: T'$
and that the composition of maps $T^1 \rightarrow T^1_{loc}
\rightarrow T'$ is surjective.

We claim that  $T^1 \rightarrow T^1_{loc}$ is in fact surjective.
To see this, it is enough to prove that $f'$ is an isomorphism
onto its image. Recall from \cite[ch.~9]{Alt} that $\dim
T^1_{loc}=9-r$ for $6\leq r \leq 8$. In the case $r=5$, the
singularity is Gorenstein  of codimension 3 so defined by
Pfaffians of a $2t+1\times 2t+1$ skew-symmetric matrix. Any
deformation of such singularity is obtained by varying its
entries. It follows that the miniversal space is smooth of
dimension $4=9-r$. Observe now, that
$$\dim T'=\dim (\ker \varphi).$$ Following
\cite[example~4.1]{Gross1} we dualize $\varphi$ and obtain the map
$$H^1(X,\Omega^1_X)\rightarrow H^0(R^1 \pi_{*}(\Omega^1_X)).$$
Since $\deg E \geq 5$
 using the Theorem on Formal Functions and the exact sequence
 $$ 0\rightarrow \Omega_X^1 \otimes
 \mathcal{I}_E^{n-1}/
\mathcal{I}_E^{n}\rightarrow \Omega_X^1 \otimes
\mathcal{O}_X/\mathcal{I}_E^{n}\rightarrow \Omega_X^1 \otimes
\mathcal{O}_X/\mathcal{I}_E^{n-1}\rightarrow 0$$ we see that
$H^0(R^1 \pi_{*}(\Omega^1_X)) \simeq H^1(E,\Omega^1_E)$, and that
 $$\varphi^{\vee}: H^1(X,\Omega^1_X)
\rightarrow H^1(E,\Omega^1_E)$$ is the restriction morphism.

Since $H^1(X,\Omega^1_X) \simeq \Pic(X)\otimes \mathbb{C}$ and
$$\dim(H^1(E,\Omega^1_E)) = \dim(\Pic(E)\otimes \mathbb{C})=10-r,$$ we
deduce from lemma \ref{pierwsze} that $$\dim (\coker \varphi
^{\vee})=9-r.$$ The claim follows.

Following the discussion on page $211$ in \cite{Gross1} and using
theorems 1.9 and 2.2 from \cite{Gross1}, we see that $$\Def (Y)=
\Spec(R/J)$$ for $R=\Lambda [[y_1,\dotsc,y_t]]$, where
$\Spec(\Lambda)$ is the base space of the miniversal deformation
of the germ $(Y',P)$, $t=\dim(\ker(T^1\rightarrow T^1_{loc}))$
moreover, $J \subset m_{\Lambda}R+ m_R^2$ is some ideal (where
$m_R=m_{\Lambda}R+(y_1,\dotsc,y_t)$).

Furthermore, $ m_{R}/( m_R^2+J) = m_{R}/ m_R^2$ (this says that
$T^1$ is isomorphic to the Zariski tangent space to $\Spec(R)$ )
and $\mbox{Supp}(R/J)=\mbox{Supp}(R)$.

 From \cite[ch.~9]{Alt} we obtain a description of $\Lambda$.
 \begin{enumerate}
\item For $r=5$ we have $\Lambda =\mathbb{C}[[x_1,x_2,x_3,x_4]]$,
\item for $r=6$ we have $\Lambda
=\mathbb{C}[[x_1,x_2,x_3]]/(x_1x_3,x_2x_3)$ then $\Spec(\Lambda)$
is a germ of a line and a plane meeting in one point, \item for
$r=7$ we have $\Lambda =\mathbb{C}[[x_1,x_2]]/(x_1^2,x_1x_2)$ and
$\Spec(\Lambda)$ is a line with a double origin, \item for $r=8$
we have two possibilities: if $E\simeq \mathbb{P}^1\times
\mathbb{P}^1$, then $\Lambda =\mathbb{C}[[x]]$; if $E\simeq
\mathbb{F}_1$, then $\Lambda =\mathbb{C}[[x]]/(x^2)$.
\end{enumerate}
It follows that in our cases
$$\Def (Y)=\Spec (R).$$
To finish the proof, consider again the local cohomology sequence
\begin{equation}\label{ciag dokladny}
H^1_E (X, \Theta_{X}) \rightarrow H^1(X ,\Theta_{X})\rightarrow
H^1(X-E ,\Theta_{X})\twoheadrightarrow T'. \end{equation} Since
$(H^1_E (X, \Theta_{X}))^{\vee}\simeq
H^0(R^2\pi_*(\Omega^1_X))=H^2(\Omega^1_E)=(H^0(\Omega^1_E))^{\vee}=0$,
we obtain that the map of germs $\Def(X)\hookrightarrow \Def(Y)$
is an embedding and its image has dimension
$t=\dim(\ker(T^1\rightarrow T'=T^1_{loc}))$ (we know that
$\Def(X)$ is smooth). Moreover, the image of the morphism
$\Def(X)\hookrightarrow \Def(Y)$ is equal set theoretically to the
principal fiber of $\Spec(R)\rightarrow \Spec(\Lambda) $. It
remains to compute the dimensions of the components of
$\Spec(\Lambda[[y_1,\dotsc,y_t]])$.
\end{proof}
\begin{rem}The above theorem is proved in the case $r=6$ in \cite[p.~762]{namikawa1} using
different methods.
\end{rem}
\begin{rem} In fact we described the miniversal Kuranishi space of
$Y$. We obtained that this space is the product of the miniversal
space of the singularity that is the cone over the contracted del
Pezzo surface with an appropriate germ of a linear space.
\end{rem}
\begin{rem} \label{remark} We compute $h^{1,2}(\mathcal{Y}_t)-h^{1,2}(X)$
in the case where $Y$ is obtained after a primitive contraction of
smooth del Pezzo surface $E$ of degree $r\leq 4$ (note that in
this case $Y$ has complete intersection singularities). In fact,
from prop.~\ref{h11} it is enough to compute the difference
$\chi(\mathcal{Y}_t)-\chi(X)$. Moreover, since we know that
$$\chi(Y)-\chi(X)=\chi(E)-1,$$ it is sufficient to find
$\chi(\mathcal{Y}_t)-\chi(Y)$.

Using known facts about Milnor numbers, we see that the last
difference between the Euler characteristics
$\chi(\mathcal{Y}_t)-\chi(Y)$ depends only on the Milnor number of
the singularity of $Y$, so depends only on the degree of $E$.
These differences will be computed using explicit examples (see
Remark \ref{hodge remark}). We obtain the following:
\begin{enumerate} \item if $r=4$, then $h^{1,2}(\mathcal{Y}_t)=
h^{1,2}(X)+7$; \item if $r=3$, then $h^{1,2}(\mathcal{Y}_t)=
h^{1,2}(X)+11$; \item if $r=2$, then $h^{1,2}(\mathcal{Y}_t)=
h^{1,2}(X)+17$; \item if $r=1$, then $h^{1,2}(\mathcal{Y}_t)=
h^{1,2}(X)+29$.\end{enumerate}
\end{rem}
Compare these results with \cite[ch.~3]{MS}.

            \section{The double cover of a Fano threefold} \label{section fano}
In this section, we obtain Calabi--Yau threefolds containing a
given del Pezzo surfaces by resolving double coverings of Fano
threefolds of index $\geq 2$ branched over singular divisors.
Recall that the index of a Fano manifold $F$ is the biggest number
$g$ such that there exists an effective divisor $E$ with
$-K_F=gE$.

In construction $1$, $2$, and $3$ the covered Fano threefold has index
$2$, $3$, and $4$ respectively. If the index is $3$ the Fano threefold is
isomorphic to a quadric in $\mathbb{P}^4$. If the index is $4$ it is
isomorphic to $\mathbb{P}^3$.

\subsection{Construction $1$}

 Let $F$ be a smooth Fano manifold with Picard number
            $\rho(F)=1$ and index $2$.
            From \cite[table~12.2]{IP}, $F$ is one of the
threefolds presented in Table $2$.

            \begin{table}[h]\label{table 1}
\caption{Fano threefolds}
            \begin{center}
            \renewcommand*{\arraystretch}{1.3}
            \begin{tabular}{|c|c|c|}\hline
            $H^3=r$&$F$&$\chi(F)$                       \\ \hline
                                                                   \hline
            $1$& a generic hypersurface of degree $6$ in
            $\mathbb{P}(1{,}1{,}1{,}2{,}3)$&$-38$              \\ \hline
            $2$&a double cover of $\mathbb{P}^3$ branched along& \\
            &a smooth
            surface of degree $4$& $-16$               \\ \hline
            $3$& a smooth cubic $X_3\subset\mathbb{P}^4$&$\:-6$ \\ \hline
            $4$& a smooth
            intersection of two quadrics $X_{2{,}2}\subset\mathbb{P}^s$&$\:\:\:\;0$                    \\ \hline
            $5$&a linear section of the Grassmannian $G(2,5)$&
\\

            & in its Pl\"{u}cker embedding &$\:\;\:\:4$
            \\ \hline
            \end{tabular}

            \end{center}
\smallskip

            \end{table}
            In each of these cases $|-K_F|=|2H|$, for some divisor $H$ on $F$. Moreover, if $r\geq 3$ the divisor $H$ is very
            ample, for $r=2$ it gives a $2:1$ morphism, and for $r=1$ it has exactly one base point (see \cite[thm.~2.4.5]{IP}).
            By the adjunction formula a general element $D'\in|H|$ is a del Pezzo surface of degree $r$.

            \begin{lem}\label{lemma4.1}

            The normal double covering $\pi\colon X'\rightarrow F$
            branched along the  sum of $D'$ and a general $G'\in|3H|$ is a singular
            Calabi--Yau variety with a transversal $A_1$ singularity.
            \end{lem}
            \begin{proof} Observe that $D'$ and $G'$ meet
            transversally. This follows from the fact that the linear system $|3H|$ is very ample.
            The singularity on the double
            cover is thus a transversal $A_1$ singularity along the
            strict transform $K$ of the intersection $D'\cap G'$ to $X'$.
            This singularity is Gorenstein, so we compute
\[K_{X'}=\pi^*(K_F)+\pi^*\left(\tfrac{D'+G'}{2}
\right)=\pi^*K_F+\pi^*(-K_F)=0\;.\qedhere\]
            \end{proof}

The blow up $\gamma:X\rightarrow X'$ along the curve $K$ is a
crepant resolution. Hence, the threefold $X$ is a smooth
Calabi--Yau threefold. Let $E$ be the exceptional divisor of
$\gamma$. Denote by $D$ and $G$ the strict transforms of $D'$ and
$G'$. Note that $D$ is isomorphic to $D'$.

            \begin{lem}

            The rank of the Picard group of $X$ is $2$ (i.e. $\rho(X)=2$).
            \end{lem}

            \begin{proof}
            This follows from \cite[Example~1]{Cynk2} because $|G'-D'|=|2H|$ is ample.
            \end{proof}
            Denote by $H^*$ the pull-back of $H$ to $X$.
            Observe that from the adjunction formula
            $D|_{D}=K_{D}$ so $H^*+D|_D=0$. This suggest the
            following proposition.

            \begin{prop}\label{constr1}

            The linear system $|G|$ is base-point-free. Moreover,
            \begin{enumerate}
            \item if $r\geq 3$, then $\varphi_{|G|}$ is a primitive contraction. Its image is a projectively normal
            threefold of degree $3r$ in
            $\mathbb{P}^{r+2}$;
            \item if $r=2$, then the morphism $\varphi_{|G|}$ is $2:1$ onto a cubic in
            $\mathbb{P}^4$;
            \item if $r=1$, then the morphism $\varphi_{|G|}$ is $3:1$ onto $\mathbb{P}^3$.
            \end{enumerate}
            \end{prop}
            \begin{proof} Since $\pi^*(H)\simeq\pi^*(D')=D+E$ and $3\pi^*(H)\simeq\pi^*(G')=G+E$ we obtain $G\in |H^*+D|=|E+3D|.$
            From lemma \ref{lemma4.1} we conclude that $G\cap
            D=\emptyset$.
            If $r \geq 2$, then $H^*|_G$ is base-point-free, hence $|G|$ is base-point-free.
            In the case $r=1$, the base point freeness follows from the fact that a generic cubic in $\mathbb{P}(1,1,1,2,3)$ does
            not pass through the base
            point of the linear system $|H|$. Note that $H^* |_D=-K_D$
            and $D|_D=K_D$ so $(H^*+D)^3= 3r$.

            We compute $h^0(\mathcal{O}_X(G))$ from the exact sequence
            $$0 \longrightarrow \mathcal{O}_X \longrightarrow \mathcal{O}_X(G)\longrightarrow N_{G|X}\longrightarrow 0.$$
            Since $\mathcal{O}_{G}(H)=K_{G}$, we obtain
            $$h^0(\mathcal{O}_X(G))=h^0(\mathcal{O}_X)+h^0(N_{G|X}) = 1 + h^0(\mathcal{O}_{G}(H)).$$ Now from
    $$O \longrightarrow \mathcal{O}_F(-2H) \longrightarrow\mathcal{O}_F(H) \longrightarrow\mathcal{O}_{G}(H) \longrightarrow O$$
    we have $h^0(\mathcal{O}_{G}(H))=h^0(\mathcal{O}_F(H)).$

            Assume that $r\geq 3$. The embeddings from table
$1$ are given by the linear system $|H|$, thus
            $$h^0(\mathcal{O}_F(H))=r+2 .$$

The morphism $\varphi_{|H^*+D|}$ is then generically $1:1$ as it
does not factorize through $\pi$. Observe moreover, that the
morphism $\varphi_{|H^*+D|}$ does not contract any curve not
contained in $D$ (a curve contained in $E$ cannot be contracted,
since the elements of $(D+|H^*|)|_E$ separate points from
different fibers of $E \rightarrow K$ and $G\cap D=\emptyset$).

To show that the image of $\varphi_{|G|}$ is projectively normal,
we follow the proof of \cite[thm.~1.4]{GaP}.
\begin{lem}   Let $T$ be a smooth projective variety.
Let $\mathcal{L}$ be a line bundle on $T$ such that for some $n$
the bundle $\mathcal{L}^{\otimes n}$ gives a birational morphism
onto a normal variety and such that the natural map
$$S^m H^0(\mathcal{L})\longrightarrow H^0(\mathcal{L}^{\otimes
m})$$ is surjective for $m \geq 1$. Then, $\mathcal{L}$ gives a
birational morphism onto a projectively normal variety.
\end{lem}
\begin{proof}[Proof of lemma]
The problem is to show the normality of the image. Consider the
following commutative diagram:
\[ \begin{CD}
             {T}   @>{\varphi_{\mathcal{L}^{\otimes n}}}>>    {\mathbb{P}(H^0(T,\mathcal{L}^{\otimes n}))} \\
             @V{\varphi_{\mathcal{L}}}VV         @V{\alpha}VV\\
             {\mathbb{P}(H^0(T,\mathcal{L}))}   @>> {v} >
             {\mathbb{P}(S^n(H^0(T,\mathcal{L}))}.
             \end{CD}\]
From the assumptions, the image of $\alpha \circ
\varphi_{\mathcal{L}^{\otimes n}}$ is normal. It remains to
observe that $v$ is the Veronese embedding.
\end{proof}

We apply the above lemma to $\mathcal{L}=\mathcal{O}_X(G)$. By
lemma \ref{pierwsze}, there is an $n$ such that
$\mathcal{O}_X(nG)$ gives a morphism with normal image. Thus, to
prove that the assumption of the lemma holds it is enough to show
that
$$H^0(nG)\otimes
H^0(G)\longrightarrow H^0((n+1)G)$$ is surjective for $n\geq 1$.
Consider the following commutative diagram:
$${H^0(nG)\otimes H^0(\mathcal{O}_X)}
\hookrightarrow {H^0(nG)\otimes H^0(G)} \twoheadrightarrow
{H^0(nG)\otimes H^0(G|_G)} \ \ \ \ $$
$$\ \  \downarrow  \ \ \ \ \ \ \ \ \ \ \ \ \ \ \ \ \ \ \ \ \ \ \ \ \ \downarrow \ \ \ \ \ \ \ \ \ \ \ \ \ \ \ \ \ \  \ \ \ \ \ \ \ \ \ \downarrow \ \ \ \ \ \ \ \ $$
$$\ \ {H^0(nG)}\ \ \hookrightarrow\ \ \ \ \
{H^0((n+1)G)}\ \ \ \ \twoheadrightarrow \ \ \ {H^0((n+1)G|_G)}.$$

To show the surjectivity of the middle vertical map, it is enough
to show the surjectivity of the right-hand-side vertical map, as
the surjectivity of the left-hand-side map is clear. Since
$\mathcal{O}_X(G)|_G=K_G$ and the restriction map
$H^0(\mathcal{O}_X(nG))\longrightarrow H^0(nK_G)$ is surjective,
it is enough to prove that
$$H^0(nK_{G})\otimes H^0(K_{G})\longrightarrow
H^0((n+1)K_{G})$$ is surjective for all $n\geq 1$. This is clear
from the fact that $G'\subset F\subset \mathbb{P}^{r+1}$ is
projectively normal for $r\geq 3$ and it is the canonical model of
$G$.

In the case $r=1$, we have $h^0(\mathcal{O}_X(G))=4$. Since
$(H^*+D)^3=3$, we conclude that $\deg \varphi_{|H^*+D|}=3$.

In the case $r=2$, we have
$h^0(\mathcal{O}_X(G))=h^0(\mathcal{O}_F(H))+1=5 $ and
$(H^*+D)^3=6$. Since $\deg( \im \varphi_{|H^*+D|})\cdot \deg
\varphi_{|H^*+D|}=6$, to prove that $\varphi_{|H^*+D|}$ is $2:1$
it is enough to observe that $G|_{G}$, defines a $2:1$ morphism
onto $X_{3}\subset \mathbb{P}^3$.
\end{proof}

  \label{ex1} \label{example} In order to describe more precisely the images of the obtained primitive contractions, we compute
  $h^0(\mathcal{O}_X(nG))$.
  From the following exact sequences
$$0 \longrightarrow \mathcal{O}_X((n-1)G) \longrightarrow \mathcal{O}_X(nG)\longrightarrow N_{nG|X}\longrightarrow 0$$
$$0 \longrightarrow \mathcal{O}_F((n-3)H) \longrightarrow\mathcal{O}_F(nH) \longrightarrow\mathcal{O}_{G}(nH) \longrightarrow 0$$
we obtain as in the proof before
$$h^o(\mathcal{O}_X(nG))=h^o(\mathcal{O}_X((n-1)G)+h^o(\mathcal{O}_F(nH))
-h^o(\mathcal{O}_F((n-3)H)).$$ Using the Riemann--Roch theorem, we
have
\begin{align*}
\chi (\mathcal{O}_F(nH))
       &=\frac{1}{12}nH(nH+2H)(2nH+2H)+\frac{1}{24}c_1c_2+\frac{1}{12}nHc_2\\
       &=\frac{1}{6}rn(n+1)(n+2)+n+1.
\end{align*}
            By the Kodaira vanishing theorem, this implies that for $n\geq 1$
            $$h^o(\mathcal{O}_X(nG))=\frac{r}{2}n^3+(\frac{r}{2}
            +3)n.$$

The graded ring $\bigoplus_{n=0}^{\infty} H^0(\mathcal{O}_X(nmG))$
is the coordinate ring of the image of $\varphi_{|mG|}$ for $m\gg
0$. Hence the scheme $\operatorname{Proj}(\bigoplus_{n=0}^{\infty}
H^0(\mathcal{O}_X(nG)))$ is isomorphic to the Calabi--Yau
threefold obtained after the primitive contraction. The Hilbert
series of $\bigoplus_{n=0}^{\infty} H^0(\mathcal{O}_X(nG))$ is the
following
$$P(t)=\frac{t^4+(r-1)t^3+rt^2+(r-1)t+1}{(1-t)^4}.$$

 For $r=1$, we obtain
$$P(t)=\frac{(1-t^6)}{(1-t)^4(1-t^2)}.$$
This is the Hilbert series of a sextic in $\mathbb{P}(1,1,1,1,2)$.
We claim that $\bigoplus_{n=0}^{\infty} H^0(\mathcal{O}_X(nG))$ is
in fact isomorphic to the graded ring of such a variety. Let us
choose generators $x,y,z,t \in H^0(\mathcal{O}_X(G))$. Computing
the dimension of $ H^0(\mathcal{O}_X(2G))$ we see that we need an
additional generator of degree $2$. It is enough to prove now that
$\varphi_{|2G|}$ gives a primitive contraction onto a projectively
normal variety in $\mathbb{P}(H^0(\mathcal{O}_X(2G)))$. To show
this, we can follow the proof of the case $i\geq 3$ of
prop.~\ref{constr1}, knowing that $2G|_G$ gives an isomorphism
onto a projectively normal surface in
$\mathbb{P}(H^0(\mathcal{O}_G(2K_G)))$.

For $r=2$, we obtain
$$P(t)=\frac{(1-t^3)(1-t^4)}{(1-t)^5(1-t^2)}.$$
This gives the Hilbert series of a complete intersection of a
quadric and a cubic in $\mathbb{P}(1,1,1,1,1,2)$. We can prove as
before that the graded ring $\bigoplus_{n=0}^{\infty}
H^0(\mathcal{O}_X(nG))$ is in fact isomorphic to the graded ring
of such a variety.

For $r=3$, we obtain
$$P(t)=\frac{(1-t^3)(1-t^3)}{(1-t)^6}.$$
This proves in fact that the image of $\varphi_{|G|}$ is the
intersection of two cubics in $\mathbb{P}^5$. Indeed, from
$h^0(\mathcal{O}_X(3G))=54$, $h^0(\mathcal{O}_X(2G))=21$, and from
the fact that the image of $\varphi_{|G|}$ is not contained in a
quadric we conclude that it is contained in two cubics without
common component. It remains to observe that the image is of
degree $9$.

For $r=4$, we obtain
$$P(t)=\frac{(1-t^2)(1-t^2)(1-t^3)}{(1-t)^7},$$
we can prove as before that the image of $\varphi_{|G|}$  is a
complete intersection of two quadrics and a cubic in
$\mathbb{P}^6$.

For $r=5$, we obtain
$$P(t)=\frac{(-t^5+5t^3-5t^2+1)(1-t^3)}{(1-t)^8}.$$
This gives a Hilbert series of the intersection in $\mathbb{P}^7$
of a cubic and a linear section of the Grassmannian $G(2,5)$ in
its Pl\"{u}cker embedding. We shall sketch the proof that the
image $\tilde{X}$ of $\varphi_{|G|}$ is such an intersection. The
image $\varphi_{|G|}(G)$ is canonically embedded in a hyperplane
section $H'\simeq \mathbb{P}^6$. It follows that the intersection
of all quadrics $Q_G\subset \mathbb{P}^6$ containing
$\varphi_{|G|}(G)$ is a smooth linear section of the Grassmannian
$G(2,5)$ in its Pl\"{u}cker embedding ($G'$ is such a section).
Since
$H^0(\mathcal{I}_{\tilde{X}}(1))=H^1(\mathcal{I}_{\tilde{X}}(1))=0$,
the natural map
$$H^0(\mathcal{I}_{\tilde{X}}(2))\rightarrow
H^0(\mathcal{I}_{\tilde{X}} \otimes \mathcal{O}_{H'}(2))$$ is an
isomorphism. Denote by $Q$ the intersection of all quadrics
containing $\tilde{X}$. We see that $\tilde{X}$ is the
intersection of $Q$ with a cubic. We shall show that $Q$ is a $4$
dimensional Gorenstein irreducible linear section of $G(2,5)$. To
prove this, it is enough to show that the fourfold $Q$ is
irreducible, Gorenstein, subcanonical, and arithmetically
Gorenstein.

Indeed, by the result of Walter \cite[thm.~0.1]{Walter} we obtain
that $Q$ is defined by Pfaffians. From the fact that $Q$ is
arithmetically Gorenstein (i.e. $\bigoplus_j H^1(I_Q(j))=0$ and
$H^i(\mathcal{O}_Q(r))=0$ for all $r$ and $1 \leq i \leq 3$) we
moreover, obtain from \cite[p.~671]{Walter} that there is a
resolution
$$0\rightarrow \mathcal{O}_{\mathbb{P}^7}(-t)\rightarrow
\mathcal{E^{\vee}}(-t)\rightarrow \mathcal{E}\rightarrow
\mathcal{O}_{\mathbb{P}^7}\rightarrow \mathcal{O}_Q ,$$ where $t$
is an integer and $\mathcal{E}=\bigoplus_{i=1}^{2p+1}
\mathcal{O}_{\mathbb{P}^7}(a_i)$ such that $a_i$ are uniquely
determined. Since a linear section of $G(2,5)$ is arithmetically
Gorenstein and restricts to the hyperplane section as $Q$, the
fourfold $Q$ is isomorphic to the linear section of the
Grassmannian $G(2,5)$ in its Pl\"{u}cker embedding.

It remains to prove that $Q$ is irreducible, Gorenstein,
subcanonical, and arithmetically Gorenstein. Since the space of
quadrics containing $\tilde{X}$ is $5$-dimensional, each component
of $Q$ is at least one-dimensional. Moreover, we know that some
hyperplane section of $Q$ is irreducible. Together, these imply
that $Q$ is irreducible.

We show now that $Q$ is normal and Gorenstein. Since $\tilde{X}
\subset Q$ is a Cartier divisor and $\tilde{X}$ is Gorenstein, we
conclude that $Q$ is Gorenstein at points of $\tilde{X}$. Next, if
$P\in Q -\tilde{X}$, then the generic codimension $2$ linear section
through $P$ cuts $Q$ along an irreducible
surface $S$ of degree $5$ in $\mathbb{P}^5$.  \\
\emph{Claim} The surface $S$ is either a cone over a smooth
elliptic
curve or a Gorenstein del Pezzo surface.\\
Indeed, from Nagata's \cite[thm.~8]{nagata} classification of
surfaces of degree $5$ in $\mathbb{P}^5$, it is enough to prove
that $S$ is neither a projection of a rational normal scroll nor a
projection of a cone over a rational normal curve. This follows
from the fact that $S$ has a hyperplane section (corresponding to
the intersection with $H'$) which is a smooth elliptic curve, but
each hyperplane section of a rational normal scroll or a cone over
a rational normal curve is rational. The claim follows.

We conclude that $Q$ is normal and Gorenstein. To see that $Q$ is
subcanonical, observe that $Q_G$ is subcanonical and hence the
generic hyperplane section is subcanonical. From the adjunction
formula, we obtain $K_{Q_G}=K_Q +Q_G|_{Q_G}$, so $K_Q|_{Q_G}=
\mathcal{O}_{Q_G}(-2)$. We use now the Lefschetz hyperplane
theorem from \cite[thm.~6]{RS} to see that there is a unique Weil divisor
on $Q$ that restricts to $Q_G$ in such a way. It follows that
$K_Q=\mathcal{O}_Q(-3)$, hence $Q$ is subcanonical.

It remains to prove that $Q$ is arithmetically Gorenstein. From
the exact sequence
$$ 0 \longrightarrow \mathcal{I}_Q(j-1) \longrightarrow \mathcal{I}_Q(j)\longrightarrow
\mathcal{I}_{Q_G}(j)\longrightarrow 0$$ and the fact that $Q_G$ is
projectively normal (i.e. $H^1(\mathcal{I}_{Q_G}(j))=0$) we obtain
that $$ H^1(\mathcal{I}_Q(j-1)) \longrightarrow H^1(
\mathcal{I}_Q(j))$$ is surjective hence $H^1(
\mathcal{I}_Q(j))=H^1( \mathcal{I}_Q)=0$.

Similarly, from $$ 0 \longrightarrow \mathcal{O}_Q(j-1)
\longrightarrow \mathcal{O}_Q(j)\longrightarrow
\mathcal{O}_{Q_G}(j)\longrightarrow 0$$ we obtain that $$
H^i(\mathcal{O}_Q(j-1)) \longrightarrow H^i( \mathcal{O}_Q(j))$$
is surjective for $1\leq i \leq 3$. From the Serre duality
$$H^i(\mathcal{O}_Q)=H^{4-i}(\mathcal{O}_Q(-3))=0.$$ Thus $Q$ is
arithmetically Gorenstein. We have proved the following theorem.
\begin{thm}\label{hilbert series}  For $r=1$, the morphism
$\varphi_{|3G|}$ is a primitive contraction. If $r=2$, then
$\varphi_{|2G|}$ is a primitive contraction, for $r=3,4,5$ the
primitive contraction is given by $\varphi_{|G|}$. Moreover,
            \item{}

\begin{enumerate}
\item if $r=1$, then $\varphi_{|3G|}(X)$ is isomorphic to a sextic
in $\mathbb{P}(1,1,1,1,2)$; \item if $r=2$, then
$\varphi_{|2G|}(X)$ is isomorphic to a complete intersection of a
cubic and a quartic in $\mathbb{P}(1,1,1,1,1,2)$; \item if $r=3$,
then $\varphi_{|G|}(X)$ is isomorphic to a complete intersection
of two cubics in $\mathbb{P}^5$; \item if $r=4$, then
$\varphi_{|G|}(X)$ is isomorphic to a complete intersection of two
quadrics and a cubic in $\mathbb{P}^6$; \item if $r=5$, then
$\varphi_{|G|}(X)$ is isomorphic to a complete intersection of a
cubic and a codimension $2$ linear section of the Grassmannian
$G(2,5)$ embedded by the Pl\"{u}cker embedding in $\mathbb{P}^9$.

\end{enumerate}
            \end{thm}
\begin{rem}\label{hodge remark}  From the geometric
descriptions contained in the previous theorem we can compute the
differences between the Euler characteristics
$\chi(\mathcal{Y}_t)-\chi(X)$ (cf. Remark \ref{remark}). Since the
values of $\chi(\mathcal{Y}_t)$ are known, it is enough to find
$\chi(X)$. We do this using the commutative diagram from
\cite{Cynk}

             \[ \begin{CD}
             {X}   @>{\gamma}>>    {X'}\\
             @VVV         @VVV\\
             {\tilde{F}}   @>> {\alpha} > {F}
             \end{CD}\]
where $\alpha$ and $\gamma$ are blowings up of $F$ and $X'$ along
$D' \cap G'$ and its strict transform to $X'$, respectively. We
need also 
the additivity of the Euler characteristic
            $$\chi(\tilde{F})=\chi(F)+\chi (G \cap D)$$
     $$\chi(X)=2\chi(\tilde{F})-\chi(D')-\chi(G').$$

     \end{rem}
 \subsection{Construction $2$}

Let $F_1$ be a Fano manifold with Picard number $1$ and Gorenstein
index greater than $2$. From tables in chapter 12 of \cite{IP} we
have $ F_1=\mathbb{P}^3$ or $F_1$ is a quadric in $\mathbb{P}^4$.

$(a)$ First, let $ F_1=\mathbb{P}^3 $. In this case
$|-K_{\mathbb{P}^3}|=|-4H|$ where $H$ is a hyperplane in
$\mathbb{P}^3$. Let us denote by $D'_i$ , for $i=1{,}2{,}3$
generic (in particular smooth) elements of $|iH|$. We have
             \begin{enumerate}   \item $D'_1=\mathbb{P}^2$,
\item
             $D'_2=\mathbb{P}^1\times\mathbb{P}^1$,
                \item $D'_3$ is the blowing up of $\mathbb{P}^2$ in
                $6$
                points.\end{enumerate}
Let us moreover, choose for $i=1{,}2{,}3$ generic elements
$G'_i\in|(8-i)H|$.

            \begin{lem}

            The normal double covers $X'_i$ of $\mathbb{P}^3$
            branched along $G'_i+D'_i$ are
            Calabi--Yau manifolds with
            transversal $A_1$ singularities along the strict transforms of $G'_i\cap D'_i$.
            \end{lem}

In the same way as in construction $1$, we resolve
$\gamma:X_i\rightarrow X'_i$ and obtain Calabi--Yau manifolds with
$\rho(X_i)=2$. As before, let $D_i$, $G_i$ be the strict
transforms of $D'_i$ and $G'_i$ to $X_i$. Observe that for each
divisor $V\in Pic(X_i)$, the restriction $V|_{D_i}$ is a multiple
of $K_{D_i}$. From lemma \ref{pierwsze}, there exists a type II
primitive contraction that contracts $D_i$. To be more precise the
following theorem holds.

\begin{thm}\label{4.3}

            The linear system $|G_i|$ is base point free and the corresponding morphism is a primitive
            contraction. Moreover,
\begin{enumerate}
\item   the image of
            $\varphi_{|G_1|}$ is a threefold of degree
            $63$ in $\mathbb{P}^{20}$,

\item   the image of
            $\varphi_{|G_2|}$ is a threefold of degree $24$ in
            $\mathbb{P}^{10}$, which is an intersection of a cone over the Veronese embedding of $\mathbb{P}^{3}$ in $\mathbb{P}^{9}$ with a generic
            sextic passing through the vertex,
\item   the  image of
            $\varphi_{|G_3|}$ is a quintic
            in $\mathbb{P}^4$.
\end{enumerate}
            \end{thm}
            \begin{proof}
Observe that $G_i \in |E+ \frac{8-i}{i}D_i|$. From the adjunction
formula
                $$D_i|_{D_i}=K_{D_i}{,}\quad E|_{D_i}=-\frac{8-i}{i}K_{D_i}.$$
We compute as before the Hilbert series of $\bigoplus_{n=0}^{\infty}
H^0(\mathcal{O}_X(nG_i))$ and obtain:
\begin{enumerate}
 \item If $i=3$, then
$$P(t)=\frac{(1-t^5)}{(1-t)^5}.$$

\item If $i=2$, then
$$P(t)=\frac{(1-t^3)(1+6t+t^2)}{(1-t)^5}.$$ This is the Hilbert series of the double Veronese embedding of a sextic in
$P(1,1,1,1,2)$.

\item If $i=1$, then
$$P(t)=\frac{t^4+17t^3+27t^2+17t+1}{(1-t)^4}.$$
The proof now follows as in example \ref{example}.\qedhere
\end{enumerate}
            \end{proof}
\begin{rem}
Denote by $\mathcal{Y}_t$ some smooth Calabi--Yau manifold from
the smoothing family of the images $Y_i$ of primitive contractions
constructed in the previous theorem. We can compute Euler
characteristics $ \chi(\mathcal{Y}_t)$ without using
thm.~\ref{4.3}. To obtain $\chi(X_i) $ we follow the method
described in remark \ref{hodge remark}. The computation is
illustrated in table $3$.

           \begin{table}[h]\label{table  2}
 \caption{illustration}
            \begin{center}

            \renewcommand*{\arraystretch}{1.3}
            \begin{tabular}{|c|c|c|c|c|c|c|c|}\hline
            $i$ & $D'_i$ & $G'_i $ & $\chi(D'_i)$ & $\chi(G'_i)$ & $\chi (G'_i \cap D'_i)$ & $\chi(X_i)$ & $\chi(\mathcal{Y}_t)$  \\ \hline
            $1$& $\mathbb{P}^2$& septic & $3$&$189$&$-28$&$-240$& -            \\ \hline
            $2$&$\mathbb{P}^1\times \mathbb{P}^1$&sextic& $4$ &$128$& $-48$&$-200$& $-204$             \\ \hline
            $3$&cubic &quintic& $9$&$55$&$-60$& $-176$&$-200$ \\ \hline

            \end{tabular}
            \end{center}

\end{table}
            We find values of $ \chi(\mathcal{Y}_t)$ using the
            results from section \ref{smoothing}.
            \end{rem}

$(b)$ In the case where $F_1$ is a quadric in $\mathbb{P}^4$, we
have $|-K_F|=|3H|$, where $H\subset \mathbb{P}^4$ is a hyperplane
section. Choose del Pezzo surfaces $D'_4\in |2H|$, $D'_5 \in |H|$
(here $D'_4$ is of degree $4$ and $D'_5$ is isomorphic to
$\mathbb{P}^1 \times \mathbb{P}^1$) and generic elements $G'_5\in
|5H| $, $G'_4\in |4H|$. We construct analogously as before
Calabi--Yau threefolds $X_4$ and $X_5$ with Picard group of rank
$2$, and divisors $G_i,D_i\subset X_i$.
\begin{rem} Using the results from section \ref{smoothing}, we see that
the threefold $\varphi_{|G_5|}(X_5)$ can be smoothed and that the
Hodge numbers of this smoothing are $h^{1,1}=1$ and $h^{1,2}=101$.
This lead us to the following theorem.
\end{rem}
\begin{thm} \label{7.1}
            The linear system $|G_i|$ is base point free and $\varphi_{|G_i|}$
            is a primitive contraction contracting $D_i$.
            \begin{enumerate}
            \item The image of $\varphi_{|G_4|}$ is
            an
            intersection of a quadric and a quartic
            in $\mathbb{P}^5$.
            \item The image of $\varphi_{|G_5|}$ is a
            threefold of degree $40$ in $\mathbb{P}^{14}$ that is a degeneration of a family of quintics embedded by the Veronese embedding in
            $\mathbb{P}^{14}$.
\end{enumerate}
\end{thm}
\begin{proof}The proof is similar to the proof of thm.~\ref{hilbert
series}.
\end{proof}
\begin{rem} We cannot make analogous constructions for other Fano threefolds.
 It follows from the tables 12.3--12.6 in \cite{IP}
 that all remaining Fano threefolds have a fibration cutting
non-canonically the embedded del Pezzo surface.
\end{rem}

\section{Construction using complete intersection}\label{section complete}

In this subsection, we show another way to construct primitive
contractions. The Calabi--Yau threefold that contains the del
Pezzo surface is obtained as a resolution of a nodal complete
intersection of hypersurfaces that contains the del Pezzo surface
in its anti-canonical embedding.

\subsection{Construction $3$}\label{constr st4 w kwin}

       Let $D'\subset \mathbb{P}^4$ be an anti-canonically embedded smooth del Pezzo
       surface of degree $4$. The surface
        $D'$ is a complete intersection of two quadrics
       $f_1{,}f_2$.

       Consider the quintic $p = f_1g_1 + f_2g_2$, where $g_1{,}g_2$ are generic cubics.
       Then, $X_1=\{p=0\}$ is a Calabi--Yau manifold with $36$ nodes at the
       points where
       $f_1=g_1=f_2=g_2=0$ (taking $f_1,g_1,f_2,g_2$ as analytic coordinates near a singular point, $p=0$ is
       an equation of an ordinary double point).

       The smooth surface $S'=\{g_1=f_2=0\}$ pas through all nodes of $X'$. The blowing up of $S'\subset X'$
       is a small resolution of $X'$.

            \begin{rem}\label{flop}
            Since $D'\subset X'$ is a smooth Weil divisor passing
            through the nodes, blowing up $D'$ produces a small
            resolution. We can obtain $X$ by flopping the {36} exceptional lines.
\end{rem}
            \begin{lem}\label{kaehler}

            The Rank of the Picard group of $X$ is $2$.
            \end{lem}
            \begin{proof} Since $X$ is a small resolution of a nodal threefold, it is
            enough to show (see \cite[thm.~2]{Cynk}) that the defect $\delta$ of
            $X'$ is $1$. We know that it is greater than $0$ since $D'$ is a Weil
            divisor.

            Recall that the defect of $X'$ is equal to
            $$ h^0(\mathcal{O}_{\mathbb{P}^4}(5) \otimes
            \mathcal{I}_{S'}) - h^0(\mathcal{O}_{\mathbb{P}^4}(5)) +\mu,$$
            where $\mu$ is the number of nodes and
            $\mathcal{I}_{S'} \subset \mathcal{O}_{\mathbb{P}^4}$
            is the ideal sheaf of these nodes. Since
            $\mathcal{I}_{S'}$ is a complete intersection of two
            quadrics and two cubics, Hilbert polynomial
            computations yield $\delta=1$.
            \end{proof}
            Let $D\subset X$ be the strict transform of $D'$ and $H$ the strict transform of the hyperplane section
            to $X$. Since $S'$ is transversal to $D'$
            we see that $D'$ is isomorphic to $D$.
            The surface $D\subset X'$ is a del Pezzo surface such that $D|_D=
            K_{D}$ (by the adjunction formula) and $H|_D=-K_{D}$ (because $D'$ is
            anti-canonically embedded). Thus, $D$ can be contracted in a primitive way.

            \begin{thm}\label{8.1}
            The linear system $|D+H|$ is base point free. Moreover, $\varphi_{|D+H|}$ is
            a primitive contraction with exceptional locus $D'$. Its image is a complete intersection
            of two cubics in $\mathbb{P}^5$.
            \end{thm}

\begin{proof} Let $S\subset X'$ be the strict transform of $S'$ and
$G'=\{f_3=g_3=0\}$. Observe that $|2H|=|S+D|$, as $f_2$ cuts $X_1$
along $S$ and $D$. Analogously $|3H|=|G+S|$, so $G\in|H+D|$. Since
$G\cap D=\emptyset$ we conclude that $|D+H|$ is base point free.
Moreover, $|D+H|$ gives a birational morphism. This follows from
the fact that $|H|$ separates points from $X \setminus D$ and that
$$(D+H)^3=K_D^2+ 3K_D^2-3K_D^2+H^3=9.$$
To compute the Hilbert series of $\bigoplus_{n=0}^{\infty}
H^0(\mathcal{O}_X(nG))$ we find
$$h^0(\mathcal{O}_{G'}(K_{G'}))=h^0(\mathcal{O}_{G'}(1)).$$ Now,
since $G'\subset \mathbb{P}^4$ is projectively normal, we conclude
as in the proof of thm.~\ref{hilbert series}.
\end{proof}
\begin{rem}  We can perform an analogous construction for del Pezzo surfaces of degree $3$.
Let $D'\subset \mathbb{P}^4$ be the intersection of a generic
cubic $c$ and a hyperplane $l$. We find as before a quintic $X'$
with $24$ node defined by $q=c g_1+lg_2$, where $g_1$ is a generic
quadric and $g_2$ a generic quartic. We construct as before a
Calabi--Yau threefold $X$ with Picard group of rank $2$ and
divisors $D,H\subset X$.
\begin{thm}
The linear system $|D+H|$ is base point free. Moreover,
$\varphi_{|D+H|}$ is the primitive contraction and its image is a
complete intersection of a quadric and a quartic in
$\mathbb{P}^5$.
 \end{thm}

\end{rem}
\subsection{Construction $4$}

 The Grassmannian $G(2{,}5)=G \subset\mathbb{P}^9$ (and so
also its generic linear section) is the zero locus of $4\times4$
Pfaffians of a generic skew-symmetric $5\times5$ matrix with
linear entries.

We fix an anti-canonically embedded del Pezzo surface $D' \subset
\mathbb{P}^5$ of degree $5$ defined by Pfaffians
$p_1,p_2,p_3,p_4,p_5$, where each $p_i$ is obtained by deleting
the $i$-th row and column from a skew-symmetric matrix $M$.

            \begin{lem}
            The intersection of two generic cubics $c_1,c_2$ from the ideal of $D'$
            is a Calabi--Yau threefold $X'\subset\mathbb{P}^5$
            with $28$ ordinary double points.
            The blowing up of $D'\subset X'$ is a small resolution.
            Flopping the $28$ exceptional lines we obtain a Calabi Yau
            threefold $X$, with two-dimensional K\"{a}hler cone, containing a del Pezzo surface of degree $5$ in a
            ``primitive''
            way.
            \end{lem}
            \begin{proof} Denote by $l_1,\dotsc,l_5,t_1,\dotsc,t_5$ linear
            forms such that
$$c_1=l_1p_1-l_2p_2+l_3p_3-l_4p_4+l_5p_5, $$ $$c_2=t_1p_1-t_2p_2+t_3p_3-t_4p_4+t_5p_5.$$
Consider the following  $7\times7$ skew-symmetric matrix $N$.

$$N=\left( \begin{array}{c|c}
\begin{array}{cc}
0&*\\
 &0
\end{array}
&
\begin{array}{ccccc}
t_1&t_2&t_3&t_4&t_5\\
l_1&l_2&l_3&l_4&l_5\\
\end{array}\\
\hline
\\
\begin{array}{cc}
&\\
&\\
&\\
&\\
&
\end{array}
& \text{\Huge{\emph{M}}}

\end{array}\right).$$

\vskip15pt
 Using the fact that the Pfaffian can be expanded along any
row, we see that $c_1=P_1$ and $c_2=P_2$, where $P_i$ is the
Pfaffian obtained by deleting the $i$-th row and column from $N$.
To see that the singularities of $X'$ are nodes, we compute the
jacobian matrix of $(c_1,c_2)$. The entries of this matrix define
an ideal $J$ with support at the singularities of $X'$. Let $I$ be
the ideal defined by $p_1,\dotsc,p_5,P_3,P_4,\dotsc,P_7$. By
simple computation in Singular we prove that $J\subset I$ and that
$I$ defines a 0-dimensional scheme of degree $28$. Moreover,
making computations in finite characteristic we see that $I$ is a
radical ideal. It remains to observe that
$c_2(\mathcal{I}_{D'}/\mathcal{I}^2_{D'}(3))=28$ (these arguments
are made more precise in \cite[thm.~1]{K}).

To show that $\rho(X)=2$, we first compute
$$ h^0(\mathcal{O}_{\mathbb{P}^5}(3) \otimes \mathcal{I}) - h^0(\mathcal{O}_{\mathbb{P}^5}(3)) +28=1.$$
Since we know that $I$ is the ideal of the set of nodes we
conclude arguing as in the proof of theorem 1 from \cite{Cynk}
that $\rho(X)\leq 2$ (the assumptions are not satisfied, so we
obtain only an inequality). Now since $X$ is a small resolution
and $X'$ is not $\mathbb{Q}$-factorial, $\rho (X)\geq
2$.\end{proof}The strict transform $D$ of $D'$ on $X$ is a del
Pezzo surface that can be contracted in a primitive way. Denote by
$H\subset X $ the strict transform of a hyperplane section of
$\mathbb{P}^5$.
\begin{rem} Using thm.~\ref{hodge}, we compute that the Hodge numbers of the
smoothing of the threefold obtained after the primitive
contraction are $h^{1,1}=1$ and $h^{1,2}=50$. This suggests that
the threefold in $\mathbb{P}^6$ is defined by the $6\times 6$
Pfaffian of a generic $7\times 7$ skew-symmetric matrix.
\end{rem}
\begin{thm}
 The morphism $\varphi_{|H+D|}$ is a primitive contraction. Moreover, its
image is a threefold of degree $14$ in $\mathbb{P}^7$ defined by
the $6\times 6$ Pfaffians of a $7\times 7$ skew-symmetric matrix .
\end{thm}
\begin{proof} Since the $6\times 6$ Pfaffians of the matrix
$N$ define a projectively normal surface  $G'\subset X'\subset
\mathbb{P}^6$ such that its strict transform in $X$ is an element
of $|H+D| $, the first part follows as in theorem \ref{8.1}. To
prove that the image is defined by Pfaffians, we use the theorem
of Walter \cite[thm~0.1]{Walter} (the image have Gorenstein
singularities). We know that this image has degree $14$ so we
conclude as in \cite{Tonoli} Section 2.
\end{proof}

\begin{rem} In an analogous way, we can construct primitive
contractions of del Pezzo surfaces of degrees $3$ and $4$. Let
$D_1'\subset \mathbb{P}^5$ be the intersection of a cubic with two
hyperplanes (resp.~$D_2'\subset \mathbb{P}^5$ be the intersection
of two quadrics $q_1,q_2$ with a hyperplane). If we choose two
generic cubics $c_1,c_2\in H^0(\mathcal{I}_{D_1'}(3))$
(resp.~$c_1,c_2\in H^0(\mathcal{I}_{D_2'}(3))$), then the
equations $c_1=c_2=0$ define a Calabi--Yau threefold $X_1'$ with
$12$ (resp.~$X'_2$ with $20$) ordinary double points. We construct
as before threefolds $X_i$ with Picard group of rank $2$ and
divisors $D_i,H_i\subset X_i$.
\begin{thm}
The linear system $|D_i+H_i|$ is base point free. Moreover,
$\varphi_{|D_i+H_i|}$ is a primitive contraction and
\begin{enumerate}
\item if $i=1$, its image is a complete intersection of two
quadrics and a cubic in $\mathbb{P}^6$; \item if $i=2$, the image
is defined in $\mathbb{P}^6$ by $4\times 4$ Pfaffians of a
$5\times 5$ skew-symmetric matrix with one row and one column of
quadrics and with the remaining entries being linear forms.
\end{enumerate}
\end{thm}
\end{rem}

\vskip10pt
\begin{minipage}{4.5cm} Grzegorz Kapustka \\
Jagiellonian University\\
ul. Reymonta 4\\
30-059 Krak\'{o}w\\
Poland\\

\end{minipage}
\hfill
\begin{minipage}{4.5cm}
 Micha\l\ Kapustka\\
 Jagiellonian University\\
 ul. Reymonta 4\\
 30-059  Krak\'{o}w\\
 Poland\\

\end{minipage}\\
email:\\
Grzegorz.Kapustka@im.uj.edu.pl\\
Michal.Kapustka@im.uj.edu.pl
\end{document}